\documentclass[12pt]{article}

\usepackage{amsmath,amsfonts,amssymb}
\usepackage{bbm}

\newtheorem{defn}{Definition}[section]

\newtheorem{lem}[defn]{Lemma}
\newtheorem{prop}[defn]{Proposition}

\newtheorem{rem}[defn]{Remark}
\newtheorem{exam}[defn]{Example}
\newenvironment{proof}{{\bf Proof }}{{\vskip 0.1cm \hfill$\Box$}}

\begin{document} 

\noindent
{\Large\bf Conservativeness of non-symmetric diffusion processes generated by perturbed divergence forms}

\bigskip
\noindent
{\bf Masayoshi Takeda}
{\footnote{Supported in part by Grant-in-Aid for Scientific Research (No.22340024 (B)), Japan Society for the Promotion of Science} 
{\bf and Gerald Trutnau}
{\footnote{Supported by the Research Settlement Fund for New Faculty and the research project  \rq\rq Advanced Research and Education of Financial Mathematics\lq\lq \
at Seoul National University.}
}}\\

\noindent
{\small{\bf Abstract.}  Let $E\subset \mathbb{R}^d$, $d\ge 2$ be an unbounded domain that is either open or closed. 
If it is closed we assume that the boundary is locally the boundary of an extension domain. 
We present conservativeness criteria for (possibly reflected) diffusions with state space $E$ and  generator $L$, which in the interior of  $E$ is given 
in the following suggestive form
\begin{eqnarray*}
Lf=\frac12 \sum_{i,j=1}^{d}\partial_j(a_{ij}\partial_i f)+\sum_{i=1}^{d}B_i\,\partial_i f.
\end{eqnarray*}
Here the diffusion matrix $(a_{ij})$ is allowed to be non-symmetric, is merely assumed to consist of measurable functions, 
and satisfies locally a strict ellipticity condition. 
$B=(B_1,...,B_d)$ is a divergence free vector field that satisfies some sector condition. 
Our main tool is a recently extended forward and backward martingale decomposition, which reduces to the well-known 
Lyons-Zheng decomposition in the symmetric case.
}

\section{Introduction}
Conservativeness criteria for diffusion operators with smooth coefficients are well-known to exist. For instance, there is Feller's and Hasminskii's test for explosions (see \cite{mckean}). 
However, these tests can not be applied to diffusions with non-regular coefficients. The aim of this work is to develop conservativeness criteria 
for a general class of such kind of diffusions that may also have an additional reflection term.\\ 
Let $E\subset \mathbb{R}^d$, $d\ge 2$ be an unbounded domain that is either open or closed.  In this article we derive general conservativeness 
criteria for divergence form operators. These divergence form operators (with or without boundary conditions) may be written 
informally as
\begin{eqnarray}\label{divform}
Lf=\frac12 \sum_{i,j=1}^{d}\partial_j(a_{ij}\partial_i f),
\end{eqnarray}
where $A=(a_{ij})_{1\le i,j\le d}$,  is a matrix of locally integrable functions $a_{ij}:E\to \mathbb{R}$ with symmetric part 
$\widetilde a_{ij}:=\frac12(a_{ij}+a_{ji})$ and anti-symmetric part $\stackrel{\vee}{a}_{ij}:=\frac12(a_{ij}-a_{ji})$.\\
Our techniques apply to diffusion matrices $A$ with fairly general symmetric part, and non-symmetric part that can be written as 
$\stackrel{\vee}{a}_{ij}=b_{ij}+c_{ij}$, where $b_{ij}=-b_{ji}$ is locally bounded, $c_{ij}=-c_{ji}$, and $\partial_j c_{ij}\in L^2_{loc}$. 
To illustrate the last let us make the simplifying assumption that everything is regular enough 
and that boundary terms disappear. Then, for smooth functions $f,g$, through integration by parts we get
\begin{eqnarray*}
-\int_E  L f\,g\,dx &=& \frac12 \sum_{i,j}\int_E (\widetilde a_{ij}+b_{ij})\partial_i f\partial_j g\,dx 
-\frac12 \sum_{i}\int_E \left (\sum_{j}\partial_j c_{ij}\right)\partial_i f \,g\,dx.\\
\end{eqnarray*}
In particular, $B=(B_1,...,B_d)$, where $B_i:=\sum_{j}\partial_j c_{ij}$ has weak divergence zero w.r.t. $dx$. 
Indeed (under the simplifying assumption) we can check
\begin{eqnarray*}
-\int_E \langle B,\nabla f \rangle  \,dx =-\int_E \sum_{i,j} \partial_{j} c_{ij} \,\partial_{i}f \,  dx=
\int_E \left (\sum_{i>j}  c_{ij}+\sum_{i<j}  c_{ij}\right )\partial_{ij} f \,  dx=0, 
\end{eqnarray*}
since  $c_{ij}=-c_{ji}$.
Summarizing, the main point is that 
\begin{eqnarray}\label{representation}
-\int_E  L f\,g\,dx ={\cal{E}}^{A}(f,g)-\int_E  \langle B, \nabla f\rangle g\,dx=:{\cal{E}}^{A,B}(f,g),
\end{eqnarray}
so that divergence form operators as in (\ref{divform}) may be well studied as perturbations ${\cal{E}}^{A,B}$ of 
non-symmetric Dirichlet forms ${\cal{E}}^{A}$ with divergence free vector fields $B$. 
In fact, our viewpoint is to analyze the non-symmetric bilinear form ${\cal{E}}^{A,B}$ with general divergence 
free vector fields $B$, not necessarily related to $A$. 
Studying the bilinear form has the advantage that the assumptions on $A$, $B$, as 
specified in the section  \ref{two} below, can be fairly weak. 
Although it lacks symmetry, another fortunate feature of ${\cal{E}}^{A,B}$ is, that its co-form 
$$
{\cal{E}}^{A^*,-B}(f,g):={\cal{E}}^{A^*}(f,g)+\int_E  \langle B, \nabla f\rangle g\,dx,
$$
where ${\cal{E}}^{A^*}(f,g):={\cal{E}}^{A}(g,f)$, and $A^*=(a_{ji})_{1\le i,j\le d}$ denotes the transposed matrix,  
has the same structural properties as ${\cal{E}}^{A,B}$. This is the main replacement tool for the missing symmetry and will be used to obtain locally the Lyons-Zheng decomposition for ${\cal{E}}^{A,B}$, 
hence also for ${\cal{E}}^{A^*,-B}$ (cf. further below).\\
In order to keep the technical subtleties low, in this article we restrict our attention to sectorial $B$ (cf. condition {\bf (P4)} in section \ref{two} below). 
In this case  ${\cal{E}}^{A,B}$ as well as its co-form are non-symmetric (sectorial) Dirichlet forms.   
However, it is also possible to consider the so far most general bilinear form approach where ${\cal{E}}^{A,B}$ as well as its co-form are 
merely generalized Dirichlet forms (but see Remark \ref{restrictiveGDF}). \\
In order to obtain conservativeness criteria for (\ref{divform}) or more generally for non-symmetric bilinear 
forms as in (\ref{representation}) we adapt the probabilistic method that has been developed in \cite{Ta89}. 
Our main technical tool is the non-sectorial Lyons-Zheng decomposition which has been 
recently extended in the framework of generalized Dirichlet forms (see \cite{Tr9}, and (\ref{1}) below).  \\
Let us briefly describe 
our method, the organization of the paper, and some main results:\\   
In section \ref{two} we develop the framework. We state the exact conditions on the domain $E$, and the coefficients given by the diffusion matrix $A=(a_{ij})_{1\le i,j\le d}$, and the vector field $B$. \\
In section \ref{three} we make use of a one-point compactification procedure which is made possible by the help of the local Poincar\'e inequality. 
More precisely, for any open and bounded subset $U$ of $E$ with positive measure we can consider the part forms of ${\cal{E}}^{A,B}$ and of its co-form on $U$. Since 
both forms have the same structural properties, e.g. satisfy the (local) Poincar\'e inequality on $U$,  both can be made simultaneously conservative on the one-point compactification of $U$ w.r.t. to some cemetery $\overline{\Delta}$ (see Proposition \ref{onepoint}, and Remark \ref{dual}). Consequently, the corresponding conservative diffusion processes admit both a Lyons-Zheng decomposition as in (\ref{1}) below. 
One of the authors in \cite{Ta89} used the Lyons-Zheng decomposition for the reflected diffusion process on $\bar{U}$. We'd like to emphasize that a one-point compactification procedure 
is an easier way to construct a conservative diffusion process for non-symmetric Dirichlet forms. 
The local Poincar\'e inequality nearly always holds. In particular it holds if the boundary of $E$ is locally Lipschitz (see Remark \ref{LP}(ii)).\\
In section \ref{sectioncons} we adapt the martingale method developed in \cite{Ta89} to the non-symmetric case.
For this we need to find a nice function $\rho$ that is positive, continuous, locally in the domain of ${\cal{E}}^{A,B}$, and satisfies $\lim_{x\to \Delta}\rho(x)=\infty$ where $\Delta$ 
is the cemetery of the diffusion associated to ${\cal{E}}^{A,B}$. Then, we apply locally on the one-point compactification of the open sets  $U=U_r^{\rho}=\{x\in E:\rho(x)< r\}$ the Lyons-Zheng decomposition for 
$\rho$ and seek for growth conditions on the coefficients that imply non-explosion. The martingale parts, which also contain the symmetric part of drift, are handled as in \cite{Ta89}. 
The remaining anti-symmetric part of drift has to be treated separately. This is developed by introducing the conditions {{\bf (S1)-(S3)}} in section \ref{sectioncons}. \\
In section \ref{five} we present examples. We illustrate the first one which is without reflection. 
If $E=\mathbb{R}^d$, $\widetilde A=(\widetilde a_{ij})_{1\le i,j\le d}$, and $A$, $B$ satisfy the conditions {\bf (P1)-\bf (P4)}
of section \ref{two} with
$$
\partial_j \hspace*{-0.1cm}\stackrel{\vee}{a}_{ij} \in L^1_{loc}(\mathbb{R}^d,dx),
$$
for all $i,j$, and
$$
\beta:= B +\frac12\left (\sum_{j=1}^d \partial_j \hspace*{-0.1cm}\stackrel{\vee}{a}_{1j},...,  \sum_{j=1}^d \partial_j \hspace*{-0.1cm}\stackrel{\vee}{a}_{dj}\right ),
$$
then the diffusion associated to ${\cal{E}}^{A,B}$ (and hence the corresponding generator) is conservative if there exists a constant $M>0$ such that for a.e. $x\in \mathbb{R}^d$\\
\begin{eqnarray*}
\frac{\langle \widetilde A(x)x, x \rangle}{|x|^2+1}+\langle\beta(x),x\rangle^+\le M(|x|^2+1)\left (\log(|x|^2+1)+1\right ).\\ \nonumber
\end{eqnarray*}
This growth condition may be best compared to \cite[Proposition 1.10]{St3} where a similar result is obtained by a completely different method in the framework of generalized Dirichlet forms. 
The "quantitative" difference is that we only remain in the framework of sectorial forms and do not consider densities to the Lebesgue measure in this article.  
However, we can allow for more general, merely measurable coefficients, and may even allow for an unbounded set of singular points of $\beta$, i.e. points around which $\beta$ is unbounded (cf. condition {\bf (S1)} in section \ref{sectioncons}, Remark \ref{sing}, and (\ref{cons(b)}) in section \ref{five}).\\
Subsequently, examples corresponding to variable oblique reflection in a wedge and in the $3$-dimensional upper half-space are studied. For the details we refer to section \ref{five}.

\section{Framework}\label{two}
Let $E$ be a locally compact separable metric space, and $m$ (the reference measure)
be a $\sigma$-finite measure on $E$. \\
For a non-symmetric Dirichlet form $({\cal{E}},{\cal{ F}})$ on $L^2(E,m)$ (see \cite{mr}, \cite{o88}) and $f,g\in {\cal{ F}}$ define its symmetric part by
$$
\widetilde {\cal{E}}(f,g):=\frac12 ({\cal{E}}(f,g)+{\cal{E}}(g,f)),
$$ 
its anti-symmetric part by
$$
{\stackrel{\vee}{{\cal{E}}}}(f,g):=\frac12 ({\cal{E}}(f,g)-{\cal{E}}(g,f)),
$$
and its co-form by
$$
\hat{\cal{E}}(f,g):=\widetilde {\cal{E}}(f,g)-{\stackrel{\vee}{{\cal{E}}}}(f,g).
$$
Suppose that ${\cal{E}}$ as well as its co-form $\hat{\cal{E}}$ are associated to a conservative diffusion. Then (see \cite{Tr9}) the {\it (extended) Lyons-Zheng decomposition} holds for any $u\in {\cal{ F}}$:
\begin{eqnarray}\label{1}
\widetilde u(X_t)-\widetilde u(X_0) & = & \frac12 M_t^{[u]}-\frac12\left \{ \hat{M}_T^{[u]}(r_T)-
\hat{M}_{T-t}^{[u]}(r_T)\right \}\nonumber \\ 
&&+\frac12\left \{N_t^{[u]}-\hat{N}_t^{[u]}\right \}; \ \ 0 \le t\le T, \ \ 
P_m\mbox{{\it -a.e.}}\\ \nonumber
\end{eqnarray}
Here $r_T$ is the time reversal operator, $M_t^{[u]}$  
(resp. $\hat{M}_t^{[u]}$) is the MAF of finite energy, and 
$N_t^{[u]}$ (resp. $\hat{N}_t^{[u]}$) is the CAF of zero energy appearing in the Fukushima decomposition (see \cite[\mbox{Theorem 4.5}]{Tr1})
corresponding to ${\cal{E}}$ (resp. $\hat{\cal{E}}$). In fact, it has been shown in \cite{Tr9} that (\ref{1}) holds in the framework of generalized Dirichlet forms, but in that case one has to 
be careful and  to precise the range of functions for which (\ref{1}) is valid (see \cite{Tr9}). \\
\begin{rem}
(i) If ${\cal{E}}$ is symmetric then $M^{[u]}=\hat{M}^{[u]}$, and $N^{[u]}=\hat{N}^{[u]}$. Thus decomposition (\ref{1}) extends the well-known Lyons-Zheng decomposition for symmetric Dirichlet forms obtained in \cite{LZ} to the non-sectorial case (for details we refer to \cite{Tr9}).\\
(ii) For a general viewpoint on decompositions of type (\ref{1}) but with non-martingale part being of bounded variation, as well as to some related stochastic calculus, we refer to \cite{RuVAWo}.\\
\end{rem}
If 
$$
{\cal F}_{loc}:=\{f:E\to \mathbb{R}\,|\, \forall G\subset E, G \mbox{ relatively compact open, }\exists g\in{\cal F} \mbox{ with } f=g \ m\mbox{-a.e. on } G \},
$$
then ($\ref{1}$) extends to $u\in {\cal F}_{loc}$ by arguments which are similar to those used in the symmetric case (cf. \cite{Tr2}).
\subsection{Non-symmetric Dirichlet forms}
{\bf (a) Strongly local forms without first order perturbation}\\ \\
Let $E\subset \mathbb{R}^d$, $d\ge 2$, with interior $E^0\not=\emptyset$, boundary $\partial E$, and $\overline{E}:=E^0\cup \partial E$. We always assume that $E^0$ is connected. 
Let $dx$ denote the Lebesgue measure on $\mathbb{R}^d$. Assume either {\bf (O)} or {\bf (C)}:\bigskip
\begin{itemize}
\item[{\bf (O)}] $E$ be open in $\mathbb{R}^d$ and $D:=C_0^{1}(E)$, i.e. $D$ is the set of continuously differentiable functions with compact suppport in $E$.\medskip
\item[{\bf (C)}] $E$ is closed in $\mathbb{R}^d$, $dx(\partial E)=0$, and $D:=\{f:E\to \mathbb{R}: \exists u\in C_0^{1}(\mathbb{R}^{d}) \mbox{ with } u=f \mbox{ on } E\}$.	
\end{itemize}
\bigskip
Let $A=(a_{ij})_{1\le i,j\le d}$,  where $a_{ij}:E\to \mathbb{R}$ with symmetric part $\widetilde a_{ij}:=\frac12(a_{ij}+a_{ji})$ and anti-symmetric part $\stackrel{\vee}{a}_{ij}:=\frac12(a_{ij}-a_{ji})$. 
Denote by $|\cdot|=\langle\cdot,\cdot\rangle^{\frac12}$ the euclidean norm on $\mathbb{R}^d$. By abuse of notation, we also use $|\cdot|$ for the absolute value in $\mathbb{R}$, i.e. in case $d=1$
Consider the following assumptions on $a_{ij}$:\bigskip
\begin{itemize}
\item[{\bf (P0)}] $a_{ij}\in L^1_{loc}(E,dx),\  1\le i,j\le d$.\medskip
\item[{\bf (P1)}] For any compact set $K\subset E$, there is a positive constant $\delta(K)$ such that 
$$
\sum_{i,j=1}^d \widetilde a_{ij}\xi_i\xi_j\ge \delta(K)|\xi|^2 
$$
$dx\mbox{-a.e. on }\ K \ \ \forall \xi=(\xi_1,...,\xi_d)\in \mathbb{R}^d$. 
\medskip
\item[{\bf (P2)}] There is a constant $M$ such that for any compact set $K\subset E$ we have 
$$
|\stackrel{\vee}{a}_{ij}|\le  \delta(K) M 
$$
$dx$-a.e. on  $K$,  $1\le i,j\le d$.
\end{itemize}
\bigskip
Then under either assumption {\bf (O)} or {\bf (C)} the bilinear form
\begin{eqnarray}\label{aijform}
{\cal{E}}^{A}(f,g):=\frac12\sum_{i,j=1}^d\int_E a_{ij}\,\partial_i f \,\partial_j g\,dx ; \ \ \ f,g\in D
\end{eqnarray}
is closable in $L^2(E,dx)$. In fact this follows from \cite{roeschmu} in case of {\bf (O)}, and in case of {\bf (C)} the arguments are similar to those used in \cite{roeschmu} (cf. e.g. \cite[Lemma 1.1]{Tr3}). \\
Let $(\cdot,\cdot)$ be  the inner product in $L^2(E,dx)$, and ${\cal{E}}^{A}_{\alpha}(f,g):={\cal{E}}^{A}(f,g)+\alpha(f,g)$, $f,g\in D$, $\alpha>0$.
Denote the  closure of $D$ w.r.t. $\|\cdot\|:={\cal{E}}^{A}_1(\cdot,\cdot)^{\frac12}$ by $D({\cal{E}}^{A})$. \\
For the anti-symmetric part we obtain 
by assumptions {\bf (P1)}, {\bf (P2)}
$$
|{\stackrel{\vee}{{\cal{E}}^{A}}}(f,g)|\le M d\sqrt{{\cal{E}}^{A}(f,f)}\sqrt{{\cal{E}}^{A}(g,g)}, \ \ f,g\in D
$$
Therefore ${\cal{E}}^{A}$ satisfies the {\it strong sector condition} 
$$
|{\cal{E}}^{A}(f,g)|\le M(1+d)\sqrt{{\cal{E}}^{A}(f,f)}\sqrt{{\cal{E}}^{A}(g,g)}, \ \ f,g\in D
$$
and can uniquely be extended outside the diagonal to $D({\cal{E}}^{A})$. In particular $({\cal{E}}^{A},D({\cal{E}}^{A}))$ is a non-symmetric regular Dirichlet form. \\ \\
Let $A^{*}:=(a_{ji})_{1\le i,j\le d}$ be the transposed matrix of $A$. Of course, $A$ satisfies {\bf (P0)}, {\bf (P1)}, {\bf (P2)}, iff $A^{*}$ does. 
Thus $({\cal{E}}^{A^*},D)$ is also closable in $L^2(E,dx)$. The closure 
$({\cal{E}}^{A^*},D({\cal{E}}^{A^*}))$, is just the co-form of $({\cal{E}}^{A},D({\cal{E}}^{A}))$. Of course it is also a regular non-symmetric Dirichlet form that satisfies the strong sector condition. 
We have
$$
D({\cal{E}}^{A})=D({\cal{E}}^{A^*}).
$$
{\bf (b) Strongly local forms perturbed by divergence free vector fields}\\ \\
We keep the notions and assumptions of {\bf (a)}. Let $B:E\to \mathbb{R}^d$ be a locally integrable vector field with the following property
\begin{itemize}
\item[{\bf (P3)}] For all $u\in D$ we have 
$$
\int_E \langle B,\nabla u \rangle  dx =0. 
$$
\end{itemize}
Let $({\cal{E}}^{A},D({\cal{E}}^{A}))$ be the Dirichlet form of {\bf (a)}.
Since closability depends only on the symmetric part we obtain under assumption {\bf (P3)} that 
$$
{\cal{E}}^{A,B}(u,v):={\cal{E}}^{A}(u,v)-\int_E \langle B,\nabla u \rangle v\,dx ; \ \ \ u,v\in D
$$
is also closable in $L^2(E,dx)$. Suppose that $B$ is additionally sectorial in the sense that there is a positive constant $C$ such that 
\begin{itemize}
\item[{\bf (P4)}] For all $u\in D$ we have 
$$
\left |\int_E \langle B,\nabla u \rangle v dx\right |\le C\sqrt{{\cal{E}}^{A}(u,u)}\sqrt{{\cal{E}}^{A}_1(v,v)}. 
$$
\end{itemize}
Note that on the right hand side there is the term ${\cal{E}}^{A}(u,u)$ and not ${\cal{E}}^{A}_1(u,u)$. This leads to a strong sector condition in $u$ and a weak 
sector condition in $v$, and vice versa. In particular, by {\bf (P4)} the statement of Proposition \ref{recurrent} below follows easily and holds also for the dual process. \\
Under  {\bf (P4)} the closure $({\cal{E}}^{A,B},{\cal F})$ of $({\cal{E}}^{A,B},D)$ in $L^2(E,dx)$ is a non-symmetric regular Dirichlet form. In particular ${\cal F}=D({\cal{E}}^{A})$,
\begin{eqnarray}\label{diagonal}
{\cal{E}}^{A,B}(u,u)={\cal{E}}^{A}(u,u)=\widetilde {\cal{E}}^{A,B}(u,u)=\widetilde {\cal{E}}^{A}(u,u)={\cal{E}}^{\widetilde A}(u,u), \ \ \forall u\in {\cal F},
\end{eqnarray}
and we have a "strong/weak" sector condition
\begin{eqnarray}\label{sector}
\left | {\cal{E}}^{A,B}(u,v)\right  |\le \left (M(1+d)+C\right )\sqrt{{\cal{E}}^{A}(u,u)}\sqrt{{\cal{E}}^{A}_1(v,v)}, \ \ \forall u,v\in {\cal F}.
\end{eqnarray}
\begin{rem}\label{restrictiveGDF}
Condition ${\bf (P4)}$ is a restrictive assumption. It can be removed in the framework of generalized Dirichlet forms. However, in what follows we will make use of quite  
many concepts which very likely hold, but still are not developed in the theory of generalized Dirichlet forms. For instance the concept of the part (and part process) of  a generalized Dirichlet form 
is not developed. Therefore we remain in the sectorial framework and assume ${\bf (P4)}$. Roughly, one can say that generalized Dirichlet form theory allows $B=(B_1,...,B_d)$, with $B_i\in L^2_{loc}(E,dx)$ (see e.g. \cite{St3}, \cite{Tr3},
whereas the sectorial framework only allows 
$B_i\in L^d_{loc}(E,dx)$, $1\le i\le d$ (see Example \ref{sobolev}).
\end{rem}
\begin{exam}\label{sobolev}
Let $d\ge 3$, $E=\mathbb{R}^d$, $K_n:=\{|x|\le n\}$, and  $B=(B_1,...,B_d)$, $B_i\in L^d_{loc}(\mathbb{R}^d,dx)$, $1\le i\le d$ with the property that for any $i$ there exists a 
constant $C_i>0$ such that for any $n\ge 1$
\begin{eqnarray}\label{sobolev1}
\min\left \{\|B_i^2\|_{\infty,K_n},\|B_i\|_{d,K_n}\right \}\le C_i\delta(K_n),
\end{eqnarray}
where $\delta(K_n)$ is the constant appearing in ${\bf (P1)}$, and $\|B_i^2\|_{\infty,K_n}$ denotes 
the essential $\sup$-norm of $B_i^2$ on $K_n$, and $\|B_i\|_{d,K_n}$ denotes the $L^d$-norm of $B_i$ on $K_n$.\\
Let $u,v\in C_0^1(\mathbb{R}^d)$, and let $n$ be such that $K_n\supset K^0_n\supset supp(u)\cup supp(v)$. Let $C_s$ be the Sobolev constant in $\mathbb{R}^d$. Then 
$$
\left |\int_{\mathbb{R}^d}B_i\partial_i u\, v \,dx\right |\le \max \left  \{\sqrt{C_i}, C_s\cdot C_i \right \}\sqrt{{\cal{E}}^{A}(u,u)}\sqrt{{\cal{E}}^{A}_1(v,v)}.
$$
Indeed if $\|B_i^2\|_{\infty,K_n}\le \|B_i\|_{d,K_n}$, then by (\ref{sobolev1}), Cauchy-Schwarz, and ${\bf (P1)}$  
\begin{eqnarray*}
\left |\int_{\mathbb{R}^d}B_i\partial_i u\, v\, dx\right | &\le &\sqrt{\int_{K_n}C_i\delta(K_n)(\partial_i u)^2  dx}\sqrt{\int_{\mathbb{R}^d}v^2 dx}\\
&\le & \sqrt{C_i}\sqrt{{\cal{E}}^{A}(u,u)}\sqrt{{\cal{E}}^{A}_1(v,v)}.\\
\end{eqnarray*}
If $\|B_i^2\|_{\infty,K_n}\ge \|B_i\|_{d,K_n}$ then similarly, applying additionally H\"older's inequality with $\frac{1}{d/2}+\frac{1}{p^*/2}=1$, we get
\begin{eqnarray*}
\left |\int_{\mathbb{R}^d} B_i\partial_i u \, v\, dx\right | &\le & \sqrt{\int_{K_n}  B_i^2\, v^2\, dx}\frac{1}{\sqrt{\delta(K_n)}}\sqrt{{\cal{E}}^{A}(u,u)}\\
&\le & \|B_i\|_{d,K_n}\|v\|_{p^*,K_n^0}\frac{1}{\sqrt{\delta(K_n)}}\sqrt{{\cal{E}}^{A}(u,u)}.\\
\end{eqnarray*}
By Sobolev's inequality and ${\bf (P1)}$ 
$$
\|v\|_{p^*,K_n^0}\le C_s \| \nabla v\|_{2,K_n^0}\le C_s\frac{1}{\sqrt{\delta(K_n)}} \sqrt{{\cal{E}}^{A}(v,v)}.
$$ 
Thus applying (\ref{sobolev1}) we obtain ${\bf (P4)}$.
\end{exam}
\section{One-point compactification and local Poincar\'e inequality}\label{three}
Let $U\subset \mathbb{R}^d$ be an arbitrary relatively compact open set with $dx(E\cap U)>0$. Let
$$
D_U:=C_0^1(U\cap E)\ \ \mbox{in case of {\bf (O)}},
$$
and
$$
D_U:=\{f:E\cap U\to \mathbb{R}: \exists u\in C_0^{1}(U) \mbox{ with } u=f \mbox{ on } E\cap U\}\ \ \mbox{in case of {\bf (C)}}.
$$
$f\in D_U$ has compact support in $E\cap U$ and one may extend $f$ to $E$ by letting $f\equiv 0$ on $E\cap U^c$. With this trivial extension we obviously obtain $D_U\subset D$.  
Since $({\cal{E}}^{A,B},D)$ is closable in $L^2(E,dx)$, we then 
have that $({\cal{E}}^{A,B},D_U)$ is closable in $L^2(E\cap U,dx)$. Denote the closure by $({\cal{E}}^{A,B},{\cal F}_U)$.
We assume that the {\it Poincar\'e inequality} holds for $({\cal{E}}^{A,B},{\cal F}_U)$, i.e. \\
\begin{itemize}
\item[{\bf (LP)}] There is a constant $c>0$, depending only on the dimension $d$, and $U$, such that for all $u\in D_U$ we have 
$$
\int_{U\cap E}u^2dx\le c\,{\cal{E}}^{A,B}(u,u)\left (=c\,{\cal{E}}^{A}(u,u)\right ). 
$$
\end{itemize}
\centerline{}
{\bf (LP)} is an abbreviation for {\it local Poincar\'e inequality}. 
\centerline{}
\begin{rem}\label{LP}
(i) The reason for assuming {\bf (LP)} is twofold. First, it allows to construct a conservative diffusion on the one-point compactification of $E\cap U$, 
whose dual is also a conservative diffusion (see Proposition \ref{onepoint}). We can hence locally apply 
 the decomposition (\ref{1}) by choosing appropriate open sets $U=U_r^{\rho}$ and let $r\to\infty$ (see section \ref{sectioncons}). 
 The second reason is that we avoid technical difficulties 
which occur when one considers the reflected diffusion on the closure of $E\cap U$. 
For the diffusion on the one-point compactification of $E\cap U$, no additional boundary terms appear on $\partial (E\cap U)\setminus \partial E$.  \\ \\
(ii) In case of {\bf (O)}, {\bf (LP)} always holds. In case of {\bf (C)}, if $E^0\cap U$ is a (bounded) extension domain, then {\bf (P1)} implies {\bf (LP)}. 
Indeed, suppose that {\bf (LP)} does not hold. Then there exists a sequence $(u_n)_{n\in\mathbb{N}} \subset {\mathcal F}_U$ 
such that $\int_{U\cap E}u_n^2dx=1$ and ${\mathcal E}^A(u_n,u_n)\to 0$ as $n\to\infty $.   
Let $G$ be a relatively compact open set with $G\supset E\cap U$.  
Since ${\mathcal F}_U$ is included in  $H^1(E^0\cap U)$, there exists $(\widetilde{u}_n)_{n\in \mathbb{N}}\subset H^1_0(G)$ 
such that $\widetilde{u}_n=u_n$ on $E\cap U$ and $\sup_n\|\widetilde{u}_n\|_{H_0^1(G)}<\infty $. 
Since the embedding $H_0^1(G)\hookrightarrow L^2(G,dx)$  is compact, there exists a subsequence $(\widetilde{u}_k)_{k\in \mathbb{N}} $ 
of $(\widetilde{u}_n)_{n\in\mathbb{N}} $ and $u\in L^2(G,dx)$ 
such that $\int_G(\widetilde{u}_k-u)^2dx \to 0$, in particular $\int_{E\cap U}
({u}_k-u)^2dx \to 0$. Note that $(u_k)_{k\in \mathbb{N}} $ is an ${\mathcal E}^A$-Cauchy 
sequence because ${\mathcal E}^A(u_n,u_n)\to 0$, and thus    
$u$ belongs to $ {\mathcal F}_U$ and ${\mathcal E}^A(u,u)=0$. Hence $u=0$ by the transience of 
 $({\mathcal E}^{\widetilde A}, {\mathcal F}_U)$, which is contradictory to $\int_{E\cap U}u^2dx=1$.
\end{rem}
\centerline{}
\centerline{}
We adjoin an extra point $\overline{\Delta}$ to $E\cap U$ and consider the one-point compactification as topology on $(E\cap U)_{\overline{\Delta}}:=(E\cap U)\cup \{\overline{\Delta}\}$. Any function $f$ on $E\cap U$ is considered as a function on 
$(E\cap U)_{\overline{\Delta}}$ by setting $f(\overline{\Delta})=0$.
We extend $dx$ to $(E\cap U)_{\overline{\Delta}}$ by setting $dx(\{\overline{\Delta}\}):=0$. Define
$$
\overline{D}_U:=\{u=u_0+k:(E\cap U)_{\overline{\Delta}}\to\mathbb{R} ; u_0\in D_U, k\in\mathbb{R} \mbox{ is a constant}\}.
$$
By {\bf (LP)} the representation of a function in $\overline{D}_U$ is unique, i.e. if $u_0+k= \widetilde u_0+\widetilde k$ then $u_0= \widetilde u_0$ in $L^2((E\cap U)_{\overline{\Delta}},dx)$. Thus\\
\begin{eqnarray}\label{compactification}
\overline{{\cal{E}}}^{A,B}(u,v):={\cal{E}}^{A,B}(u_0,v_0);\ \ \ u,v\in \overline{D}_U
\end{eqnarray}
is well-defined.\\
\begin{prop}\label{onepoint}
$(\overline{{\cal{E}}}^{A,B},\overline{D}_U)$ is closable in $L^2((E\cap U)_{\overline{\Delta}},dx)$. The closure
$(\overline{{\cal{E}}}^{A,B},\overline{{\cal{ F}}}_U)$ is a  strongly 
local conservative non-symmetric regular Dirichlet form. The associated process can be decomposed as in (\ref{1}). $({\cal{E}}^{A,B},{\cal F}_U)$ 
is the part Dirichlet form on $E\cap U$ of $({\cal{E}}^{A,B},{\cal F})$, as well as of $(\overline{{\cal{E}}}^{A,B},\overline{{\cal{ F}}}_U)$. Moreover 
$u\in \overline{{\cal{ F}}}_U$, iff $u=u_0+k$ for some $u_0\in {\cal{ F}}_U$ and $k\in \mathbb{R}$.
\end{prop}
\centerline{}
\begin{rem}\label{dual}
Since $A^*$ (resp. $-B$) satisfy the same assumptions as $A$ (resp. $B$) the corresponding statements of Proposition \ref{onepoint} also hold for the Dirichlet forms  
$({\cal{E}}^{A^*,-B},{\cal F}_U)$,  $({\cal{E}}^{A^*,-B},{\cal F})$, and $(\overline{{\cal{E}}}^{A^*,-B},\overline{{\cal{ F}}}_U)$.
\end{rem}
\centerline{}
\begin{proof}(of Proposition \ref{onepoint})
Let $(u_n:=u_0^n+k_n)_{n\in \mathbb{N}}\subset \overline{D}_U$ be $\overline{{\cal{E}}}^{A,B}$-Cauchy such that $u_n\to 0$ in $L^2((E\cap U)_{\overline{\Delta}},dx)$. Then by {\bf (LP)} and definition of $\overline{{\cal{E}}}^{A,B}$
$$
\int_{U\cap E}|u^n_0-u^m_0|^2dx\le c\,{\cal{E}}^{A,B}(u^n_0-u^m_0,u^n_0-u^m_0)=c\,\overline{{\cal{E}}}^{A,B}(u^n-u^m,u^n-u^m).
$$
Therefore $(u^n_0)_{n\in \mathbb{N}}$ converges in $L^2(E\cap U,dx)$ and is ${\cal{E}}^{A,B}$-Cauchy. Since $(u_n)_{n\in \mathbb{N}}$ converges to zero,  $(u^n_0)_{n\in \mathbb{N}}$ must converge to a constant, say $k\in \mathbb{R}$. Since ${\cal{E}}^{A,B}$ is closable 
we obtain
$$
\lim_{n\to\infty}\overline{{\cal{E}}}^{A,B}(u_n,u_n)=\lim_{n\to\infty}{\cal{E}}^{A,B}(u_0^n,u_0^n)={\cal{E}}^{A,B}(k,k)=0.
$$
Thus $(\overline{{\cal{E}}}^{A,B},\overline{D}_U)$ is closable in $L^2((E\cap U)_{\overline{\Delta}},dx)$. By (\ref{sector}), and (\ref{compactification}), $(\overline{{\cal{E}}}^{A,B},\overline{D}_U)$ satisfies a sector condition. 
Therefore $(\overline{{\cal{E}}}^{A,B},\overline{{\cal{ F}}}_U)$ is a coercive closed form. In particular $(v_n:=v_0^n+l_n)_{n\in \mathbb{N}}\subset \overline{D}_U$ is $\overline{{\cal{E}}}_1^{A,B}$-Cauchy, iff $(v_0^n)_{n\in \mathbb{N}}\subset D_U$ is ${\cal{E}}_1^{A,B}$-Cauchy. Therefore the last statement of the Proposition follows.\\
For any $\varepsilon>0$ let $\varphi_{\varepsilon}:\mathbb{R}\to [-\varepsilon,1+\varepsilon]$,   $\varphi_{\varepsilon}\in C_b^1(\mathbb{R})$,  $\varphi_{\varepsilon}(t)=t$, $t\in[0,1]$, $\varphi_{\varepsilon}'\in [0,1]$,    $\varphi_{\varepsilon}(t)=1+ \varepsilon $, $t\ge 1+2\varepsilon$,    $\varphi_{\varepsilon}(t)=-\varepsilon $, $t\le -2\varepsilon $. 
Let $v=v_0+l\in \overline{D}_U$. We have
\begin{eqnarray*}
\varphi_{\varepsilon}(v) & = & \underbrace{\varphi_{\varepsilon}(v_0+l)- \varphi_{\varepsilon}(l)}_{\in D_U}+\underbrace{\varphi_{\varepsilon}(l)}_{\in\mathbb{R}}\in \overline{D}_U.\\
\end{eqnarray*}
The fact that $(\overline{{\cal{E}}}^{A,B},\overline{{\cal{ F}}}_U)$ is a Dirichlet form can now be derived using \cite[\mbox{I. Proposition 4.10}]{mr}
\end{proof}
\section{Conservativeness criteria}\label{sectioncons}
Let $\mathbb{M}=((X_t)_{t\ge 0},(P_x)_{x\in E\cup\{\Delta\}})$ be the diffusion associated to $({\cal{E}}^{A,B},{\cal F})$.
In the following we want to find conservativeness criteria for $\mathbb{M}$ in the cases\\
\begin{itemize}
\item[{\bf (AO)}] The conditions {\bf (O)}, {\bf (P0)}, {\bf (P1)}, {\bf (P2)}, {\bf (P3)}, {\bf (P4)}, {\bf (LP)} hold and $E\subset\mathbb{R}^d$ is open and unbounded,
\end{itemize}
and
\begin{itemize}
\item[{\bf (AC)}] The conditions {\bf (C)}, {\bf (P0)}, {\bf (P1)}, {\bf (P2)}, {\bf (P3)}, {\bf (P4)}, {\bf (LP)} hold and $E\subset\mathbb{R}^d$ is closed and unbounded.
\end{itemize}
\bigskip
Suppose that {\bf (AO)}, or {\bf (AC)} holds. \\ \\
Let $(G_{\alpha})_{\alpha>0}$ be the resolvent of $({\cal{E}}^{A,B},{\cal F})$, and $({\widehat G}_{\alpha})_{\alpha>0}$ be the co-resolvent. 
Just as in the symmetric case (see \cite[Theorem 1.6.6]{fot}), a basic conservativeness criterium is given by:\\
\begin{lem}\label{basic}
Suppose there is $v_0\in L^2(E,dx)\cap L^1(E,dx)$, $v_0>0$ $dx$-a.e, and 
$(u_n)_{n\in \mathbb{N}}\subset {\cal F}$, $0\le u_n\le 1$, $n\in \mathbb{N}$, $u_n\uparrow 1$  as $n\to\infty$, such that
$$
\lim_{n\to\infty}{\cal{E}}^{A,B}(u_n,\widehat{G}_1 v_0)=0.
$$
Then $({\cal{E}}^{A,B},{\cal F})$ is conservative.
\end{lem}
\begin{proof}
Since
$$
0=\lim_{n\to\infty}{\cal{E}}^{A,B}(u_n,\widehat{G}_1 v_0)=\lim_{n\to\infty}\int_E (u_n-G_1 u_n)v_0 dx=\int_E (1-G_1 1)v_0 dx,
$$ 
we get $G_1 1=1$ as desired.
\end{proof}
\centerline{}
\begin{prop}\label{recurrent}
(cf. \cite[Theorem 1.5.10]{o88} in case of strong sector condition) Suppose that 
the symmetric part $({\cal{E}}^{\widetilde A},{\cal F})$  of $({\cal{E}}^{A,B},{\cal F})$ is recurrent in the sense of \cite{fot}. 
Then $({\cal{E}}^{A,B},{\cal F})$ is conservative.
\end{prop}
\begin{proof}
By \cite[Theorem 1.6.3]{fot} 
there exists $(u_n)_{n\in \mathbb{N}}\subset {\cal F}$, $0\le u_n\le 1$, $n\in \mathbb{N}$, $u_n\uparrow 1$  as $n\to\infty$, such that
$$
{\cal{E}}^{\widetilde A}(u_n,u_n)\to 0 \ \ \mbox{ as } \ n\to \infty.
$$
Let $v_0$ be as in Lemma \ref{basic}. Applying the sector condition (\ref{sector}) we get
\begin{eqnarray*}
\left | {\cal{E}}^{A,B}(u_n,\widehat{G}_1 v_0)\right  |\le \left (M(1+d)+C\right )\sqrt{{\cal{E}}^{A}(u_n,u_n)}\sqrt{{\cal{E}}^{A}_1(\widehat{G}_1 v_0,\widehat{G}_1 v_0)}.
\end{eqnarray*}
By (\ref{diagonal}) and recurrence the right hand side of the last expression tends to zero as $n\to\infty$. Therefore the conservativeness of $({\cal{E}}^{A,B},{\cal F})$ now follows from 
Lemma \ref{basic}. 
\end{proof}
\centerline{}
Lemma \ref{basic} and Proposition \ref{recurrent} are quite simple. Criteria for recurrence, and non-recurrence can e.g. be found in \cite{fot}, \cite{o88}, \cite{stu}.\\ \\ 
Fix $\rho\in C(E)^+\cap {\cal F}_{loc}$ such that $\lim_{x\to \Delta}\rho(x)=\infty$. For any $r>0$ let 
$$
U_r^{\rho}:=\{x\in E:\rho(x)< r\}.
$$ 
By Proposition \ref{onepoint} and Remark \ref{dual} on each $(U_r^{\rho})_{\overline{\Delta}}$ we have the Lyons-Zheng decomposition w.r.t. $P^{r}_{m_r}$, $m_r:=dx_{|U_r^{\rho}}$, $m_r(\{\overline{\Delta}\})=0$, for the conservative diffusion $\overline{\mathbb{M}}^r$ associated to $(\overline{{\cal{E}}}^{A,B}, \overline{{\cal{ F}}}_{U_r^{\rho}})$. \\ \\
Denote the part process of $\mathbb{M}$ (and $\overline{\mathbb{M}}^r$) on $U_r^{\rho}$ by $\mathbb{M}^r=((X_t^{0,r})_{t\ge 0},(P^{0,r}_x)_{x\in U_r^{\rho}\cup\{\overline{\Delta}\}})$. 
Let $T>0$. If $X_0\in U_R^{\rho}$, then \\
\begin{eqnarray*}
C_r&:=&\{\sup_{t\in [0, T]}\left (\rho(X_t)-\rho(X_0)\right )\ge r\}
\\
&=&\{\exists t_0\in [0,T] \mbox{ with }\rho(X_{t_0})=\rho(X_0)+ r \mbox{ and } \rho(X_{t})< \rho(X_0)+ r \ \forall t< t_0\},\\
\end{eqnarray*}
and thus the event $C_r$ only depends on the behaviour of $X_t$ strictly before it hits the boundary $\partial U_{R+r}^{\rho}:=\{x\in E:\rho(x)=R+r\}$ for the first time. Hence 
$$
P_x(C_r)=P^{0,R+r}_x(C_r)=P^{R+r}_x(C_r) \ \ \ \  \mbox{ for } \  dx\mbox{-a.e. } \ x\in U_R^{\rho},
$$
and in particular
\begin{eqnarray}\label{cr}
C_r=\left \{\sup_{t\in[0, T], \, t<\tau_{R+r}}\left (\rho(X_t)-\rho(X_0)\right )\ge r\right \}
\end{eqnarray}
whenever $X_0\in U_R^{\rho}$ where $\tau_{R+r}:=\inf\{t\ge 0:X_t\notin U_{R+r}^{\rho}\}$. We obtain
$$
P_{m_R}\left (C_r\right )=\int_{U_R^{\rho}}P_x\left (C_r\right )dx=\int_{U_R^{\rho}}P^{R+r}_x\left (C_r\right )dx=P^{R+r}_{m_R}\left (C_r\right ),
$$
and thus 
\begin{eqnarray*}
P_{m_R}\left (\sup_{t\in [0, T]}\left (\rho(X_t)-\rho(X_0)\right ) = \infty\right ) & = & \lim_{r\to \infty}P_{m_R}\left (C_r\right )
 =  \lim_{r\to \infty}P_{m_R}^{R+r}\left (C_r\right )\\
&  &\hspace*{-5cm} \le \liminf_{r\to \infty}P_{m_{R+r}}^{R+r}\left (\sup_{t\in[0, T], \, t<\tau_{R+r}}(\rho(X_t)-\rho(X_0))\ge r\right ).\\
\end{eqnarray*}
The function $(\rho-(R+r))\wedge 0$ belongs to the part Dirichlet space 
on $U^\rho_{R+r}$. Thus if we define $\rho_{R+r}:=((\rho-(R+r))\wedge 0)+(R+r)$, 
$\rho_{R+r}$ belongs to $\bar{\mathcal F}_{U^\rho_{R+r}}$ and 
$\rho_{R+r}=\rho$ on $U^\rho_{R+r}$. 
Thus using in particular ($\ref{1}$), and ($\ref{cr}$) we obtain
\begin{eqnarray*}
P_{m_{R+r}}^{R+r}\left (\sup_{t\in[0, T], \, t<\tau_{R+r}}(\rho(X_t)-\rho(X_0))\ge r\right ) & \le & 
P_{m_{R+r}}^{R+r}\left (\sup_{t\in[0, T], \, t<\tau_{R+r}}\frac12 M_t^{[\rho]}\ge \frac{r}{4}\right )\\
& & \hspace*{-8cm} +P_{m_{R+r}}^{R+r}\left (\sup_{t\in[0, T], \, t<\tau_{R+r}}-\frac12\hat{M}_T^{[\rho]}(r_T)\ge \frac{r}{4}\right )+
P_{m_{R+r}}^{R+r}\left (\sup_{t\in[0, T], \, t<\tau_{R+r}}\frac12\hat{M}_{T-t}^{[\rho]}(r_T)\ge \frac{r}{4}\right )\\ 
& & \hspace*{-5cm} +P_{m_{R+r}}^{R+r}\left (\sup_{t\in[0, T], \, t<\tau_{R+r}}\frac12\left \{N_t^{[\rho]}-\hat{N}_t^{[\rho]}\right \}\ge \frac{r}{4}\right )\\
& & \hspace*{-8.5cm}\le P_{m_{R+r}}^{R+r}\left (\sup_{t\in[0, T], \, t<\tau_{R+r}}M_t^{[\rho]}\ge \frac{r}{2}\right )
+P_{m_{R+r}}^{R+r}\left (\sup_{t\in[0, T], \, t<\tau_{R+r}}\hat{M}_{t}^{[\rho]}\ge \frac{r}{2}\right )
\\
&&\hspace*{-8cm} +P_{m_{R+r}}^{R+r}\left (\sup_{t\in[0, T], \, t<\tau_{R+r}}-\hat{M}_{t}^{[\rho]}\ge \frac{r}{2}\right )+P_{m_{R+r}}^{R+r}\left (\sup_{t\in[0, T], \, t<\tau_{R+r}}\frac12\left \{N_t^{[\rho]}-\hat{N}_t^{[\rho]}\right \}\ge \frac{r}{4}\right )\\
\end{eqnarray*}
The terms with martingale part can be handled as in \cite{Ta89}, or see also \cite[5.7]{fot}. Thus for the martingale parts we need to show 
\begin{eqnarray}\label{martingaledrift}
\lim_{r\to \infty}3 \,\cdot \mbox{vol}(U_{R+r}^{\rho}) \,\cdot  \mbox{Erfc}\left (\frac{r}{\sqrt{8 M^{\rho}(R+r)\cdot T}}\right )=0,
\end{eqnarray}
where $\mbox{Erfc(x)}:=\frac{2}{\sqrt{\pi}}\int_x^{\infty} e^{-x^2}dx$, and $M^{\rho}(r):=\mbox{ess.sup}\{\langle \widetilde{A}\nabla \rho, \nabla \rho\rangle(x)\,|\, x\in U^{\rho}_r\}$. 
It remains to find additional conditions that also  ensure
\begin{eqnarray}\label{drift}
\lim_{r\to \infty}P_{m_{R+r}}^{R+r}\left (\sup_{t\in[0, T], \, t<\tau_{R+r}}\frac12\left \{N_t^{[\rho]}-\hat{N}_t^{[\rho]}\right \}\ge \frac{r}{4}\right )=0.
\end{eqnarray}
\centerline{}
In general, the anti-symmetric part of drift w.r.t. $\rho$, namely 
$$
\frac12\left \{N_t^{[\rho]}-\hat{N}_t^{[\rho]}\right \},
$$ 
has an absolutely continuous, and an non-absolutely continuous part. We assume it may be written as
$$
\int_0^t\langle \beta,\nabla \rho \rangle(X_s) ds+\int_0^t f(X_s) dG_s,
$$
where (provided the $\stackrel{\vee}{a}_{ij}$ are regular enough)
$$
\langle \beta,\nabla \rho \rangle =\sum_{i=1}^d\left ( B_i +\frac12\sum_{j=1}^d \partial_j \hspace*{-0.1cm}\stackrel{\vee}{a}_{ij}\right )\partial_i \rho,
$$
$G=G^1-G^2$ is the difference of two PCAF's, and $f$ some function. \\ 
Therefore, in order to check (\ref{drift}), it is enough to verify
\begin{eqnarray}\label{drift1}
\lim_{r\to \infty}P_{m_{R+r}}^{R+r}\left (\sup_{t\in[0, T], \, t<\tau_{R+r}}\int_0^t\langle \beta,\nabla \rho \rangle(X_s) ds\ge \frac{r}{8}\right )=0,
\end{eqnarray}
and 
\begin{eqnarray}\label{drift2}
\lim_{r\to \infty}P_{m_{R+r}}^{R+r}\left (\sup_{t\in[0, T], \, t<\tau_{R+r}}\int_0^t f(X_s) dG_s\ge \frac{r}{8}\right )=0.
\end{eqnarray}
We refer below to (\ref{drift1}) as the {\bf absolutely continuous case}, and to (\ref{drift2}) as the  {\bf non-absolutely continuous  case}. \\ \\
{\bf (a) Absolutely continuous case} \\ \\
Since $\langle \beta,\nabla \rho \rangle$ may be singular (i.e. $=+\infty$) in some points but still locally integrable in a neighborhood of such 
singular points we introduce the set $S_{int}$, i.e. we suppose that there is a measurable set $S_{int}\subset E$, with:
\begin{itemize}
\item[{\bf (S1)}] For any $R>0$
$$
\lim_{r\to\infty}\frac1r\int_{U_{R+r}^{\rho}}\langle \beta,\nabla \rho \rangle^{+}\cdot \mathbb{I}_{S_{int}}(x)dx=0.
$$ 
\end{itemize}
where $\langle \beta,\nabla \rho \rangle^{+}:=\langle \beta,\nabla \rho \rangle\cdot \mathbb{I}_{\{\langle \beta,\nabla \rho \rangle\ge 0\}}$.\\ \\
On the complement of neighborhoods of singular points we assume
\begin{itemize}
\item[{\bf (S2)}] There is  a constant $c_1$ such that (a.e.) for any $r>0$
$$
\langle \beta,\nabla \rho \rangle \cdot \mathbb{I}_{U_r^{\rho}\setminus S_{int}}\le c_1 (1+r).
$$ 
\end{itemize}
Now under {\bf (S2)}, if $0<T<\frac{1}{16 c_1}$, and $r$ is large, writing $\mathbb{I}_{U_{R+r}^{\rho}}=\mathbb{I}_{U_{R+r}^{\rho}\setminus S_{int}}+\mathbb{I}_{U_{R+r}^{\rho}\cap S_{int}}$ we get
\begin{eqnarray*}
P_{m_{R+r}}^{R+r}\left (\sup_{t\in[0, T], \, t<\tau_{R+r}}\int_0^t\langle \beta,\nabla \rho \rangle(X_s) ds\ge \frac{r}{8}\right ) & \le & 
P_{m_{R+r}}^{R+r}\left (T c_1(1+R+r)\ge \frac{r}{16}\right )\\
&& \hspace*{-4cm}+P_{m_{R+r}}^{R+r}\left (\int_0^T\langle \beta,\nabla \rho \rangle^{+} \mathbb{I}_{U_{R+r}^{\rho}\cap S_{int}}(X_s) ds\ge \frac{r}{16}\right )\\
& & \hspace*{-6cm} = P_{m_{R+r}}^{R+r}\left (\int_0^T\langle \beta,\nabla \rho \rangle^{+} \mathbb{I}_{U_{R+r}^{\rho}\cap S_{int}}(X_s) ds\ge \frac{r}{16}\right ).\\
\end{eqnarray*}
For the estimate of the last term we use the Chebyshev-Markov inequality, and Fubini's Theorem:
\begin{eqnarray*}
P_{m_{R+r}}^{R+r}\left (\int_0^T\langle \beta,\nabla \rho \rangle^{+} \mathbb{I}_{U_{R+r}^{\rho}\cap S_{int}}(X_s) ds\ge \frac{r}{16}\right )
& \le  & \frac{16}{r} E_{m_{R+r}}^{R+r}\left [\int_0^T\langle \beta,\nabla \rho \rangle^{+} \mathbb{I}_{U_{R+r}^{\rho}\cap S_{int}}(X_s) ds\right ]\\
& & \hspace*{-4cm} = \frac{16}{r} \int_0^T \int_{(U_{R+r}^{\rho})_{\overline{\Delta}}} p_s\left (\langle \beta,\nabla \rho \rangle^{+} \mathbb{I}_{U_{R+r}^{\rho}\cap S_{int}}\right )(x) dxds\\
& &  \hspace*{-4cm}= \frac{16}{r} \int_0^T \int_{(U_{R+r}^{\rho})_{\overline{\Delta}}} \langle \beta,\nabla \rho \rangle^{+} \mathbb{I}_{U_{R+r}^{\rho}\cap S_{int}}(x) \hat p_s\mathbb{I}_{(U_{R+r}^{\rho})_{\overline{\Delta}}}dxds.\\
\end{eqnarray*}
By conservativeness of the co-process to $\overline{\mathbb{M}}^{R+r}$ (see Remark \ref{dual}) we have $\hat p_s\mathbb{I}_{(U_{R+r}^{\rho})_{\overline{\Delta}}}=\mathbb{I}_{(U_{R+r}^{\rho})_{\overline{\Delta}}}$ and the last term equals 
$$
\frac{16T}{r} \int_{U_{R+r}^{\rho}} \langle \beta,\nabla \rho \rangle^{+}\cdot \mathbb{I}_{S_{int}}(x) dx.
$$
Consequently, applying additionally {\bf (S1)} we get for $0<T<\frac{1}{16 c_1}$, 
\begin{eqnarray*}
\lim_{r\to \infty}P_{m_{R+r}}^{R+r}\left (\sup_{t\in[0, T], \, t<\tau_{R+r}}\int_0^t\langle \beta,\nabla \rho \rangle(X_s) ds\ge \frac{r}{8}\right ) 
& = & 0.\\
\end{eqnarray*}
\begin{exam}
Let $E=\mathbb{R}^2$, $\stackrel{\vee}{a}_{12}=-\stackrel{\vee}{a}_{21}=|x|^{\gamma}\wedge 1$, $\gamma>0$, and $B\equiv 0$. Let $\rho(x)=\log(|x|+2)$. 
Then ({\bf S1}), {\bf (S2)}, {\bf (P2)}, are satisfied with e.g. $S_{int}=\{x\in \mathbb{R}^2\,:\, |x|\le 1\}$. Suppose that 
$(a_{ij})_{1\le i,j\le 2}$ satisfies {\bf (P0)}, {\bf (P1)}. Then we are in the situation of {\bf (AO)}, and employing 
the function $\rho(x)$ (cf. \cite{Ta89}, or \cite[5.7]{fot}), the associated process is conservative, if for a.e. $x$ 
$$
\frac{\langle \widetilde A(x)x, x \rangle}{|x|^2} \le const.(|x|+2)^2\log(|x|+2).
$$
\end{exam}
\bigskip
{\bf (b) Non-absolutely continuous  case} \\ \\
Let $f, G$ be as in (\ref{drift2}). There are mainly three types of possibilities: 
\begin{itemize}
\item[(i)] $G_s$ corresponds to a local time on the boundary $U_{R+r}^{\rho}\cap \partial E$. 
In this case, if $\stackrel{\vee}{A}$, $\partial E$,  are sufficiently regular and $\eta$ denotes 
the inward normal, we have      $f=\langle \stackrel{\vee}{A}\hspace*{-0.1cm}\eta,\nabla \rho \rangle$.
\item[(ii)] $G_s$ corresponds to a reflection term inside $U_{R+r}^{\rho}$, e.g. if (at least) one of 
the $\stackrel{\vee}{a}_{ij}$ has a jump discontinuity along some $(d-1)$-dimensional hyperspace in $B_{R+r}^{\rho}$.
\item[(iii)] The distributional derivative $\frac12\partial_j \hspace*{-0.1cm}\stackrel{\vee}{a}_{ij}$ is a 
smooth measure $\mu_{ij}$ that is not absolutely continuous w.r.t. the Lebesgue measure. 
In this case $f\equiv 1$ and
$$
G_t=\sum_{i=1}^{d}\left (\sum_{j=1}^{d}\int_0^t \partial_i \rho(X_s)dH^{ij}_s\right )
$$ 
where $H^{ij}$ is a CAF uniquely related to $\mu_{ij}$.
\end{itemize}
Denote by  $\mu_{G^i}$ the smooth measure associated to $G^i$, by $V_1(\cdot)$ the corresponding $1$-potentials, and let $g_1:=f^+$, $g_2:=f^-$. 
Then, in any of the above cases (i)-(iii)  
\begin{eqnarray*}
&&\hspace*{-2cm}P_{m_{R+r}}^{R+r}\left (\sup_{t\in[0, T], \, t<\tau_{R+r}}\int_0^t f(X_s) dG_s\ge \frac{r}{8}\right ) \\ 
&& \le \sum_{i=1}^{2} P_{m_{R+r}}^{R+r}\left (\int_0^T g_i\mathbb{I}_{U_{R+r}^{\rho}}(X_s) dG_s^i\ge \frac{r}{16}\right )\\
&&\le \frac{16e^T}{r}\sum_{i=1}^{2} E_{m_{R+r}}^{R+r}[V_1(g_i \mathbb{I}_{U_{R+r}^{\rho}}\cdot G^i)]\\ 
&& = \frac{16e^T}{r}\sum_{i=1}^{2}\int_{U_{R+r}^{\rho}}g_i d\mu_{G^i}.\\ 
\end{eqnarray*}
Thus similarly to {\bf (S1)} we introduce
\begin{itemize}
\item[{\bf (S3)}] For any $R>0$
$$
\lim_{r\to\infty}\frac1r \left ( \int_{U_{R+r}^{\rho}}f^+ d\mu_{G^1}+\int_{U_{R+r}^{\rho}}f^- d\mu_{G^2}\right )=0.
$$ 
\end{itemize}
\centerline{}
\section{Applications and Examples}\label{five}
Below we show through examples how to apply the conditions {\bf (S1)-(S3)}. 
The \lq\lq optimal\rq\rq \ result depends highly on the choice 
of the function $\rho$.  In general the function $\rho$ should be chosen case by case accordingly to the explicitly given coefficients. 
However,  a straightforward and unsophisticated application of {\bf (S1)-(S3)} with a $\rho$ of logarithmic type already leads to conservativeness results for a wide range of diffusions.\\ \\ 
{\bf (a) Example without reflection}\\ \\
Assume {\bf (AO)} with $E=\mathbb{R}^d$. By Remark \ref{LP}, {\bf (O)} always implies {\bf (LP)}. Thus the main assumptions in {\bf (AO)} 
are {\bf (P0)-(P4)}. 
In order to simplify things we assume that 
$$
\partial_j \hspace*{-0.1cm}\stackrel{\vee}{a}_{ij} \in L^1_{loc}(\mathbb{R}^d,dx),
$$
for all $i,j$, i.e. the antisymmetric part of drift is of bounded variation.\\ \\
We fix  
$$
\rho(x)=\log(|x|^2+1).\\
$$
Showing that (\ref{martingaledrift}) follows if there is some constant $M_1>0$ such that for a.e. $x$\\
\begin{eqnarray}\label{sym(b)}
\frac{\langle \widetilde A(x)x, x \rangle}{|x|^2+1}\le M_1(|x|^2+1)\left (\log(|x|^2+1)+1\right ), \\ \nonumber
\end{eqnarray}
we see that the symmetric part is conservative if (\ref{sym(b)}) holds.\\
Since
$$
\nabla \rho(x)=\frac{2x}{|x|^2+1}
$$
is locally bounded we can see that for the antisymmetric part of drift 
$$
\beta= B +\frac12\left (\sum_{j=1}^d \partial_j \hspace*{-0.1cm}\stackrel{\vee}{a}_{1j},...,  \sum_{j=1}^d \partial_j \hspace*{-0.1cm}\stackrel{\vee}{a}_{dj}\right ),
$$
$\langle \beta,\nabla \rho \rangle$ is locally integrable. Clearly, if the set of singular points $S_{int}$ is bounded then {\bf (S1)} trivially holds. \\
\begin{rem}\label{sing}
We can even allow for an unbounded set of singular points. 
In this case we need to verify {\bf (S1)}, i.e. we need to verify that for any $R>0$
$$
\lim_{r\to\infty}\frac1r\int_{\left \{|x|\le \sqrt{e^{R+r} -2}\right \}}\frac{\langle \beta(x), 2x\rangle ^+}{|x|^2+1} \cdot \mathbb{I}_{S_{int}}(x)dx=0.
$$ 
Of course this result can also be refined regarding separately the singular points of each $B_i$, $\partial_j \hspace*{-0.1cm}\stackrel{\vee}{a}_{ij}$.\\
\end{rem}
On the complement of singular points we assume: There is a constant $M_2>0$ such that \\
\begin{eqnarray}\label{anti(b)}
\langle\beta(x),x\rangle \le M_2(|x|^2+1)\left (\log(|x|^2+1)+1\right )\\ \nonumber
\end{eqnarray}
for a.e. $x\in \mathbb{R}^d\setminus S_{int}$. Using (\ref{anti(b)}) one can easily see that {\bf (S2)} holds. \\ \\
Combining  (\ref{sym(b)}) and (\ref{anti(b)}) and assuming (for simplicity) that $S_{int}$ is a bounded set, we obtain: 
If there exists a constant $M>0$ such that for a.e. $x\in \mathbb{R}^d$\\
\begin{eqnarray}\label{cons(b)}
\frac{\langle \widetilde A(x)x, x \rangle}{|x|^2+1}+\langle\beta(x),x\rangle^+ \cdot \mathbb{I}_{\mathbb{R}^d\setminus S_{int}}(x)\le M(|x|^2+1)\left (\log(|x|^2+1)+1\right ),\\ \nonumber
\end{eqnarray}
then the diffusion $\mathbb{M}=((X_t)_{t\ge 0},(P_x)_{x\in \mathbb{R}^d})$ associated to $({\cal{E}}^{A,B},{\cal F})$ is conservative.\\ \\
{\bf (b) Oblique reflection}\\ \\
{\bf I. Brownian motion with variable oblique reflection in a wedge:}\\ \\
Let $p_1, p_2\in (0,1)$, and $p_1^2+p_2^2=1$. Let 
$e_r=(p_1,p_2)$, and $e_l=(-p_1,p_2)$. Define a wedge of angle 
$$
\phi=\arccos(p_2^2-p_1^2)\in (0,\pi)
$$
by
$$
E=W_{\phi}:=\{z=(x,y)\in \mathbb{R}^2\,|\, z=ae_r+be_l;\,a,b\ge 0\}.
$$
Thus we are in the situation of {\bf (C)}. The wedge has two inward normal vectors one on the left hand side, and one on the right hand side
$$
\eta_r(z)\equiv \eta_r= (-p_2,p_1),\ \ \ \ \ \mbox{and}\ \ \ \ \ \ \eta_l(z)\equiv\eta_l=(p_2,p_1).
$$
Oblique reflection on the right hand side (resp. left hand side) of $\partial W_{\phi}$ is determined by the angle between $A(z)\eta_{r}$ and $e_r$ 
(resp. the angle between $A(z)\eta_{l}$ and $e_l$). On the right hand side of the wedge the reflection angle is given by
\begin{eqnarray}\label{ranglewedge}
\theta_r(z)=\arccos \left(\frac{|\langle A(z)\eta_{r},e_r\rangle |}{| A(z)\eta_r|}\right ).
\end{eqnarray}
and on the left hand side it is given by 
\begin{eqnarray}\label{langlewedge}
\theta_l(z)=\arccos \left(\frac{|\langle A(z)\eta_{l},e_l\rangle |}{| A(z)\eta_l|}\right ).
\end{eqnarray}
The particular case $\theta(z)=\frac{\pi}{2}$ corresponds to normal reflection at the point $z$. 
Note that the reflection angle is variable.\\
We will consider the Brownian motion case with variable oblique reflection, so 
$$
\widetilde a_{ij}=\delta_{ij}, \ \ \ \ \ \ \ i,j=1,2.
$$ 
Thus {\bf (P1)} holds. Assume that we are given a vector field $B:=(B_1,B_2)$ that satisfies {\bf (P3)}, and {\bf (P4)}, and assume also that 
{\bf (P0)}, {\bf (P2)} hold for $\stackrel{\vee}{a}_{12}$. In order to obtain oblique reflected Brownian motion, we have to eliminate the absolutely continuous part of drift. Thus we have to assume that $\stackrel{\vee}{a}_{12}$ is weakly differentiable and
\begin{eqnarray}\label{obliquedivzero}
\beta=\left (B_1+\frac12 \partial_2 \hspace*{-0.1cm}\stackrel{\vee}{a}_{12},B_2-\frac12 \partial_1 \hspace*{-0.1cm}\stackrel{\vee}{a}_{12}\right )=0.
\end{eqnarray}
Right below let us give an example where all our assumptions hold.
\begin{exam}\label{exampled=2}
Define
$$
B_1(x,y):=-j'(y)f\left (k(x)+j(y)\right ),
$$
and
$$
B_2(x,y):=k'(x)f\left (k(x)+j(y)\right ),
$$
where $f,k,j\in H^{1,1}_{loc}(\mathbb{R})$, and $B_1$, $B_2$ are supposed to be bounded. Then clearly $B:=(B_1,B_2)$ satisfies {\bf (P4)}. 
If
\begin{eqnarray*}
j'\left (-\frac{p_2}{p_1}x\right )=-j'\left (\frac{p_2}{p_1}x\right )=k'(x)
\end{eqnarray*}
then one can readily check that {\bf (P3)} holds. Let 
$$
\stackrel{\vee}{a}_{12}(x,y):=g\left (k(x)+j(y) \right ),
$$
where $g$ is bounded, weakly differentiable, and $g'=f$. Thus  {\bf (P0)}, and {\bf (P2)} hold. \\
For the non-symmetric absolutely continuous part of drift, we can easily check (\ref{obliquedivzero}).
According to  (\ref{ranglewedge}), (\ref{langlewedge}) we get 
$$
\theta_r(z)=\arccos \left(\frac{\stackrel{\vee}{a}_{12}\hspace*{-0.1cm}(z)}
{\sqrt{1+\stackrel{\vee}{a}_{12}\hspace*{-0.1cm}(z)^2}}\right ),  \  
\mbox{ and } \ \
\theta_l(z)=\arccos \left(\frac{-\stackrel{\vee}{a}_{12}\hspace*{-0.1cm}(z)}
{\sqrt{1+\stackrel{\vee}{a}_{12}\hspace*{-0.1cm}(z)^2}}\right )
$$
\centerline{}
as variable reflection angles. In particular the corresponding diffusion is conservative by Remark \ref{restrictive}(ii).
\end{exam}
\centerline{}
\begin{rem}\label{restrictive}
(i) Variably oblique reflected Brownian motion has been considered under various aspects and by more probabilistic means in \cite{rog1}, \cite{rog2}. In the case of constantly oblique reflected Brownian motion we have $\stackrel{\vee}{a}_{12}\equiv c$ for some constant $c\in \mathbb{R}$, 
$\theta_r(z)\equiv \theta_r$, $\theta_l(z)\equiv \theta_l$, and $\theta_l+\theta_r=\pi$.
Constantly oblique reflected BM was studied intensively under various aspects in \cite{vw}, \cite{williams}. There, the main parameter is
$$
\alpha:=\frac{\theta_l+\theta_r-\pi}{\phi}.
$$ 
Thus $\stackrel{\vee}{a}_{12}=c$ corresponds to $\alpha=0$. Using the underlying Dirichlet form we can directly see that the origin is not 
reached since it has zero capacity. 
This is shown in \cite{vw} for more general $\alpha\le 0$ by probabilistic means. \\
(ii) By \cite[Theorem 3]{stu} the process corresponding to $\stackrel{\vee}{a}_{12} \equiv 0$, i.e. the $2$-dimensional normally reflected BM in $W_{\phi}$, 
is recurrent (this is also proved in \cite{williams}). 
Thus by Proposition \ref{recurrent} the variably oblique reflected BM (i.e. the process corresponding to general 
$\stackrel{\vee}{a}_{12}\hspace*{-0.1cm}(z)$ satisfying (\ref{obliquedivzero}) above) 
is conservative. By this, we can also directly see that the Dirichlet form defined in \cite[Theorem 6.2]{kim} is conservative if $d=2$. \\
(iii) Applying condition {\bf (S3)} we were not successful in showing conservativeness for the most general $2$-dimensional 
oblique reflected BM that we can construct using the Dirichlet form method. 
In some sense condition {\bf (S3)} is too strong (the Chebyshev-Markov inequality that is used for the derivation of {\bf (S3)} 
leads to a quite rough estimate) for certain 
drifts that are not absolutely continuous, and so
we had to make a detour using recurrence and Proposition \ref{recurrent}. 
In dimension three and higher we are more flexible and it is easier to verify {\bf (S3)}, i.e. it is easier  to find a nice $\rho$ such that 
$\langle \stackrel{\vee}{A}\hspace*{-0.1cm}\eta\,,\nabla \rho \rangle\le 0$ for many different choices of $\stackrel{\vee}{A}$
(see the following example and in particular Example \ref{s4}). In general the function $\rho$ should be chosen individually and accordingly to the 
explicitly given coefficients.\\
However, assuming $\langle \,\hspace*{-0.1cm}\stackrel{\vee}{A}\hspace*{-0.1cm}\eta\,,\nabla \rho \rangle^+\le f(r)$ on $U_r^{\rho}$ we have 
\begin{eqnarray*} 
&&\hspace*{-2cm}P_{m_{R+r}}^{R+r}\left (\sup_{t\in[0, T], \, t<\tau_{R+r}}\int_0^t \langle \,\hspace*{-0.1cm}\stackrel{\vee}{A}\hspace*{-0.1cm}\eta\,,\nabla \rho \rangle(X_s)d\ell_s\ge r\right ) \\
&& \le P_{m_{R+r}}^{R+r}\left (\int_0^T \langle \,\hspace*{-0.1cm}\stackrel{\vee}{A}\hspace*{-0.1cm}\eta\,,\nabla \rho \rangle^+(X_s)\mathbb{I}_{U_{R+r}^{\rho}}(X_s) d\ell_s\ge r\right )\\
&&\le P_{m_{R+r}}^{R+r}\left (\int_0^T \mathbb{I}_{U_{R+r}^{\rho}}(X_s)\ell_s\ge \frac{r}{f(R+r)}\right ).\\
\end{eqnarray*}
So, one would rather need an estimate of 
$$
P_{m_{R+r}}^{R+r}\left (\ell_T\ge \frac{r}{f(R+r)} \right ).
$$
But again, it is not obvious to handle this term because of the appearance of $m_{R+r}$. A pointwise statement seems be easier to handle in this case. 
For instance, if we know that a PCAF $G_t$ is finite 
$P_x$-a.s. and if $\frac{r}{f(R+r)}\to \infty$, then
$$
\lim_{r\to \infty}P_x\left (G_T\ge \frac{r}{f(R+r)} \right )=P_x\left (G_T =\infty \right )=0.
$$ 
Pointwise statements can be considered by using Fukushima's decomposition, 
or by using the  Lyons-Zheng decomposition w.r.t. $P_{m_{R+r}}$ for the martingale parts, and  w.r.t. $P_x$ for a.e. 
$x\in E$ for the drift parts. 
Thus the example shows that for non-absolutely continuous drifts, 
such as local times 
it might be reasonable to consider the method in a pointwise setting. 
\end{rem}
\centerline{}
{\bf II. Variable oblique reflection in $3$-dimensional upper half-space:}\\ \\
We consider the upper half-space
$$
E=\{(x_1,x_2,x_3)\in \mathbb{R}^3\,|\,x_3\ge 0\}.
$$
Thus we are in the situation of  {\bf (C)}. Let $A=(a_{ij})_{1\le i,j\le d}$ satisfy {\bf (P1)}, 
and $\stackrel{\vee}{a}_{ij}$ satisfying {\bf (P2)}. Let $B=(B_1, B_2, B_3)$ be a vector field satisfying {\bf (P3)}, and {\bf (P4)}. We consider the Dirichlet form given by (\ref{aijform}). 
The associated process, is a diffusion $X_t$ (up to its lifetime), and  
the reflection term appearing in Fukushima's decomposition for $\rho$ is given by
$$
\int_0^t \langle A\eta,\nabla \rho \rangle(X_s)d\ell_s
$$
where $\ell_s$ is roughly speaking the local time on $\partial E$. In dimension three oblique reflection is most suitably determined by a vector field 
$$
(F(x),\theta(x))
$$ 
that assigns to each point $x\in \partial E$ a reflection direction $F(x)\in \partial E$ 
and a reflection angle $\theta(x)\in (0,\frac{\pi}{2}]$. The particular case $(F(x),\theta(x))=(0,\frac{\pi}{2})$ corresponds to normal reflection at the point $x$. 
Reflection angle and direction are hence variable and are calculated as follows:
Let $e_1=(1,0,0),...,e_3=(0,0,1)$ be the standard basis of $\mathbb{R}^3$. The interior normal of $\partial E$ is $e_3$. The reflection angle is given by
$$
\frac{\pi}{2}-\theta[A(x)\eta(x),e_3],
$$
where 
$\theta[A(x)\eta(x),e_3]$ denotes the angle between $A(x)\eta(x)$ and the vector $e_3$. Since $\eta(x)\equiv e_3$ we obtain
\begin{eqnarray}\label{angle}
\theta(x)=\arcsin \left(\frac{|\langle A(x)e_3,e_3\rangle |}{|A(x)e_3|}\right )=\arcsin \left(\frac{1}{\sqrt{a_{13}(x)^2+a_{23}(x)^2+1}}\right ).
\end{eqnarray}
The reflection direction is given as the orthogonal projection of $A(x)\eta(x)$ on $\partial E$, thus
\begin{eqnarray}\label{direction}
F(x)=\langle A(x)e_3,e_1\rangle e_1+ \langle A(x)e_3,e_2\rangle e_2=(a_{13}(x),a_{23}(x),0). 
\end{eqnarray}
If $\stackrel{\vee}{a}_{ij}\in H^{1,1}_{loc}(E)$, integration by parts shows that the absolutely continuous antisymmetric part of drift is given by 
$$
\beta\cdot \nabla =\left (B_1+\frac12 \left (\partial_2 \hspace*{-0.1cm}\stackrel{\vee}{a}_{12}+\partial_3 \hspace*{-0.1cm}\stackrel{\vee}{a}_{13}\right ),
B_2+\frac12 \left(\partial_1 \hspace*{-0.1cm}\stackrel{\vee}{a}_{21}+\partial_3 \hspace*{-0.1cm}\stackrel{\vee}{a}_{23}\right ),
B_3+  \frac12 \left(\partial_1 \hspace*{-0.1cm}\stackrel{\vee}{a}_{31}+\partial_2 \hspace*{-0.1cm}\stackrel{\vee}{a}_{32}\right )\right)\cdot \nabla .
$$
In order to be concrete we fix   
$$
\rho(x)=\log(|x|^2+2).
$$
By Remark \ref{LP}, {\bf (LP)} holds on each $E\cap U_{r}^{\rho}$. 
As before we see that the symmetric part is conservative if (\ref{sym(b)}) holds for a.e. $x$ in $E$. We also come to the same 
conclusions for $\beta$ as in Remark \ref{sing} and (\ref{anti(b)}). Thus conservativeness for the martingale and absolutely continuous parts 
of drift is guaranteed by (\ref{cons(b)}) if the corresponding set of singular points $S_{int}$ is bounded.\\
The antisymmetric reflection part is given by 
\begin{eqnarray}\label{ref55}
\int_0^t \langle \,\hspace*{-0.1cm}\stackrel{\vee}{A}\hspace*{-0.1cm}\eta\,,\nabla \rho \rangle(X_s)d\ell_s.
\end{eqnarray}
Note that the integral is well-defined since $\stackrel{\vee}{a}_{ij}\in H^{1,1}_{loc}(E)$ (as usual we denote the trace on $\partial E$ also by $\stackrel{\vee}{a}_{ij}$). 
For the drift (\ref{ref55}) with the explicitly given $\rho$ condition {\bf (S3)} holds if
$$
\lim_{r\to\infty}\frac{1}{r}\int_{\partial E} \frac{\langle ( \,\stackrel{\vee}{\hspace*{-0.1cm}a}_{13}\hspace*{-0.1cm}(z,0),\,\,\hspace*{-0.1cm}\stackrel{\vee}{a}_{23}\hspace*{-0.1cm}(z,0)) \, , 2z \rangle ^+}{ |z|^2+2} 
\mathbb{I}_{\{|z|\le \sqrt{e^{R+r} -2}\}}dz=0
$$ 
where $z=(x_1,x_2)$. The latter is for instance satisfied, if there is a constant $M_3>0$ with
\begin{eqnarray}\label{reflec(b)}
\langle ( \,\stackrel{\vee}{\hspace*{-0.1cm}a}_{13}\hspace*{-0.1cm}(z,0),\,\,\stackrel{\vee}{\hspace*{-0.1cm}a}_{23}\hspace*{-0.1cm}(z,0)) \, , z \rangle ^+\le M_3 f'\hspace*{-0.1cm}\left (\log(|z|^2+2)\right ) 
\end{eqnarray}
for a.e. $z\in \mathbb{R}^2$, where $f$ is some differentiable function with $\lim_{t\to\infty}f(t)t^{-1}=0$.\\ \\
As conclusion we have that if the set of singular points is bounded then (\ref{cons(b)}) and (\ref{reflec(b)}) lead to the conservativeness of
the oblique reflected diffusion associated to $({\cal{E}}^{A},{\cal F})$. If the set of singular points is unbounded, as before we have to verify the analogon of Remark \ref{sing}.\\
As an example, (\ref{reflec(b)}) is satisfied if for $dz$-a.e. $z$ outside some compact set we have
\begin{eqnarray*}
\stackrel{\vee}{a}_{13}\hspace*{-0.1cm}(z,0)x_1+\stackrel{\vee}{a}_{23}\hspace*{-0.1cm}(z,0)x_2\le 0,
\end{eqnarray*}
so for instance, if
\begin{eqnarray}\label{criterium00}
\stackrel{\vee}{a}_{13}\hspace*{-0.1cm}(z,0)=x_2 f(z),\ \ \ \  \mbox{and}\ \ \ \ \stackrel{\vee}{a}_{23}\hspace*{-0.1cm}(z,0)=-x_1 f(z),
\end{eqnarray}
for some function $f\in H^{1,1}_{loc}(\mathbb{R}^2)$ with growth conditions given accordingly to {\bf (P2)}. 
\begin{rem}
(\ref{criterium00}) corresponds to an antisymmetric reflection direction that is perpendicular to the radial direction.  
In order to obtain conservativeness results for more general antisymmetric reflection directions one has to choose $\rho$ more carefully (see Example \ref{s4} below). 
\end{rem}
Finally, we apply our result to \\ \\
\centerline{{\it Brownian motion with variable oblique reflection,}}
\centerline{} 
For this we first have to put 
$$
\widetilde a_{ij}=\delta_{ij}, \ \ \ \ \ \ \ i,j=1,2,3.
$$ 
The absolutely continuous non-symmetric part of drift has to disappear, so $\beta\equiv 0$ is our condition. One can compensate the 
absolutely continuous part of drift 
produced by $\stackrel{\vee}{A}$ with the vector field $B\not=0$ as we did before in a two dimensional wedge to obtain the most general oblique reflection. 
However, let us proceed in a slightly less general, but more simple way (cf. \cite{ag}, \cite{kim}). For this, we let $B=0,\stackrel{\vee}{a}_{12}=0$, and 
$\stackrel{\vee}{a}_{13}$, $\stackrel{\vee}{a}_{23}$ only depend on $z=(x_1,x_2)$. Thus  $\stackrel{\vee}{a}_{13}$, $\stackrel{\vee}{a}_{23}$ are constant in the variable $x_3$. If
\begin{eqnarray}\label{cond}
\partial_1 \hspace*{-0.1cm}\stackrel{\vee}{a}_{13}+\partial_2 \hspace*{-0.1cm}\stackrel{\vee}{a}_{23}=0,  
\end{eqnarray}
then $\beta =0$ and we obtain BM with variable oblique reflection according to (\ref{angle}) and (\ref{direction}). In the next example 
we show how to derive conservativeness for a concrete choice of $\stackrel{\vee}{a}_{13}$,  $\stackrel{\vee}{a}_{23}$, using condition {\bf (S3)}.
\begin{exam}\label{s4}
A concrete example satisfying (\ref{cond}) is e.g. given by 
\begin{align}\label{cond2}
\stackrel{\vee}{a}_{13}(x_1,x_2,x_3)  & = -l'(x_2)g\big (k(x_1)+l(x_2)\big ), \nonumber\\  
  \stackrel{\vee}{a}_{23}(x_1,x_2,x_3)   & = \ k'(x_1)g\big (k(x_1)+l(x_2)\big ), 
\end{align}
where $g\in C^1_b(\mathbb{R})$, and $k, l\in C^1(\mathbb{R})$ are positive, $k', l'\in C_b(\mathbb{R})$, and $k(t),l(t)\to\infty$ as $|t|\to\infty$. 
Choosing 
$$
\rho(x_1,x_2,x_3)=k(x_1)+l(x_2)+m(x_3).
$$
with $m\in C^1(\mathbb{R})$ positive, $m'\in C_b(\mathbb{R})$, and $m(t)\to\infty$ as $|t|\to\infty$, we obtain
$$
\langle \stackrel{\vee}{A}\hspace*{-0.1cm}\eta\,,\nabla \rho \rangle=\stackrel{\vee}{a}_{13}\hspace*{-0.1cm}(z,0)\partial_1 \rho(z,0)+\stackrel{\vee}{a}_{23}\hspace*{-0.1cm}(z,0)\partial_2 \rho(z,0)\le 0
$$
Then {\bf (S3)}  is satisfied with $f^+=\langle \stackrel{\vee}{A}\hspace*{-0.1cm}\eta\,,\nabla \rho \rangle^+=0$, $G^1=\ell$, $G^2\equiv 0$.
Since $|\nabla \rho|^2\le const.$ it follows that $M^{\rho}(R+r)\le const.$ 
If $k(x_1)+l(x_2)+m(x_3)\ge const. \sqrt{\log(|x|^{\varepsilon}+2)}$ for some $\varepsilon>0$, then 
$vol(U^{\rho}_{R+r})\le e^{\frac{1}{\varepsilon}\left (\frac{R+r}{const.}\right )^2}$. Thus choosing $T$ small enough in (\ref{martingaledrift}) 
conservativeness also follows for the Brownian motion part. Therefore the diffusion corresponding to the Dirichlet form (\ref{aijform}) with

\[
A =
\left( {\begin{array}{ccc}
 1 & 0 & \stackrel{\vee}{a}_{13} \\
 0 & 1 & \stackrel{\vee}{a}_{23} \\
 -\stackrel{\vee}{a}_{13} & -\stackrel{\vee}{a}_{23} & 1
 \end{array} } \right)
\]
where $\stackrel{\vee}{a}_{13}$ and $\stackrel{\vee}{a}_{23}$ are as in (\ref{cond2}) is conservative.
\end{exam}

Mathematical Institute, Tohoku University, Aoba, Sendai, 980-8578, Japan. \\
E-mail: takeda@math.tohoku.ac.jp\\ \\ 
Seoul National University, Department of Mathematical Sciences and Research Institute of Mathematics,
San56-1 Shinrim-dong Kwanak-gu, Seoul 151-747, South Korea,\\
E-mail: trutnau@snu.ac.kr
\end{document}